\documentclass[12pt]{article}
\usepackage{epic,eepic}
\usepackage[dvips]{epsfig}
\usepackage{color}

\usepackage{amsmath, amscd, amssymb}
\usepackage{longtable,textcomp}
\usepackage{cite}
\usepackage{enumerate}
\usepackage{authblk}
\usepackage{enumitem}

\def\draw #1 by #2 (#3){
  \vbox to #2{
    \hrule width #1 height 0pt depth 0pt
    \vfill
    \special{picture #3}
    }
  }
\def\scaleddraw #1 by #2 (#3 scaled #4){{
  \dimen0=#1 \dimen1=#2
  \divide\dimen0 by 1000 \multiply\dimen0 by #4
  \divide\dimen1 by 1000 \multiply\dimen1 by #4
  \draw \dimen0 by \dimen1 (#3 scaled #4)}
  }

\begin{document}

\renewcommand{\labelenumi}{\theenumi}
\newcommand{\qed}{\mbox{\raisebox{0.7ex}{\fbox{}}}}
\newcommand{\ep}{\epsilon}
\newcommand{\la}{\lambda}
\newcommand{\G}{\Gamma}
\newcommand{\AL}{\mathbb{A}}
\newtheorem{theorem}{Theorem}
\newtheorem{example}{Example}
\newtheorem{conjecture}[theorem]{Conjecture}
\newtheorem{remark}{Remark}
\newtheorem{problem}[theorem]{Problem}
\newtheorem{defin}[theorem]{Definition}
\newtheorem{definition}[theorem]{Definition}
\newtheorem{lemma}[theorem]{Lemma}
\newtheorem{corollary}[theorem]{Corollary}
\newtheorem{nt}{Note}
\newtheorem{proposition}[theorem]{Proposition}
\renewcommand{\thent}{}
\newenvironment{pf}{\medskip\noindent{\textbf{Proof}:  \hspace*{-.4cm}}\enspace}{\hfill \qed \medskip \newline}
\newenvironment{spf}{\medskip\noindent{\textbf{Sketch of proof}:  \hspace*{-.4cm}}\enspace}{\hfill \qed \medskip \newline}
\newenvironment{pft2}{\medskip\noindent{\textbf{Proof of Theorem~\ref{2}}:  \hspace*{-.4cm}}\enspace}{\hfill \qed \medskip \newline}
\newenvironment{defn}{\begin{defin}\em}{\end{defin}}{\vspace{-0.5cm}}
\newenvironment{lem}{\begin{lemma}\em}{\end{lemma}}{\vspace{-0.5cm}}
\newenvironment{cor}{\begin{corollary}\em}{\end{corollary}}{\vspace{-0.5cm}}
\newenvironment{thm}{\begin{theorem} \em}{\end{theorem}}{\vspace{-0.5cm}}
\newenvironment{pbm}{\begin{problem} \em}{\end{problem}}{\vspace{-0.5cm}}
\newenvironment{note}{\begin{nt} \em}{\end{nt}}{\vspace{-0.5cm}}
\newenvironment{exa}{\begin{example} \em}{\end{example}}{\vspace{-0.5cm}}
\newenvironment{rem}{\begin{remark} \em}{\end{remark}}{\vspace{-0.5cm}}
\newenvironment{pro}{\begin{proposition} \em}{\end{proposition}}{\vspace{-0.5cm}}

\def\diag{\text{diag}}
\def\al{\alpha}
\def\Ga{\Gamma}
\def\be{\beta}
\def\la{\lambda}
\def\ge{\geq}
\def\le{\leq}

\setlength{\unitlength}{12pt}
\newcommand{\comb}[2]{\mbox{$\left(\!\!\begin{array}{c}
            {#1} \\[-0.5ex] {#2} \end{array}\!\!\right)$}}
\renewcommand{\labelenumi}{(\theenumi)}
\renewcommand{\b}{\beta}
\newcounter{myfig}
\newcounter{mytab}
\def\mod{\hbox{\rm mod }}
\def\scaleddraw #1 by #2 (#3 scaled #4){{
  \dimen0=#1 \dimen1=#2
  \divide\dimen0 by 1000 \multiply\dimen0 by #4
  \divide\dimen1 by 1000 \multiply\dimen1 by #4
  \draw \dimen0 by \dimen1 (#3 scaled #4)}
  }
\newcommand{\Aut}{\mbox{\rm Aut}}
\newcommand{\w}{\omega}
\def\r{\rho}
\newcommand{\DbF}{D \times^{\phi} F}
\newcommand{\autF}{{\tiny\Aut{\scriptscriptstyle(\!F\!)}}}
\def\Cay{\mbox{\rm Cay}}
\def\a{\alpha}
\newcommand{\C}[1]{\mathcal #1}
\newcommand{\B}[1]{\mathbb #1}
\newcommand{\F}[1]{\mathfrak #1}
\title{On the characterization of geometric distance-regular graphs}

     \author[a]{Chenhui Lv}
     \author[a,b]{Jack H. Koolen\footnote{J.H. Koolen is the corresponding author.}}

    	\affil[a]{\footnotesize{School of Mathematical Sciences, University of Science and Technology of China, Hefei, 230026, People's Republic of China}}
	\affil[b]{\footnotesize{CAS Wu Wen-Tsun Key Laboratory of Mathematics, University of Science and Technology of China, Hefei, 230026, People's Republic of China}}

\date{\today}

\maketitle
\pagestyle{plain}
\newcommand\blfootnote[1]{%
		\begingroup
		\renewcommand\thefootnote{}\footnote{#1}%
		\addtocounter{footnote}{-1}%
		\endgroup}
	\blfootnote{E-mail addresses:   {\tt lch1994@mail.ustc.edu.cn} (C. Lv), {\tt koolen@ustc.edu.cn} (J.H. Koolen).}

\thanks{Dedicated to Professor Paul Terwilliger on the occasion of his 70th birthday.}

\begin{abstract}
In 2010, Koolen and Bang proposed the following conjecture:
For a fixed integer $m \geq 2$, any geometric distance-regular graph with smallest eigenvalue $-m$, diameter $D \geq 3$ and $c_2 \geq 2$ is either a Johnson graph, a Grassmann graph, a Hamming graph, a bilinear forms graph, or the number of vertices is bounded above by a function of $m$.
In this paper, we obtain some partial results towards this conjecture.


\end{abstract}

\section{Introduction}
All graphs considered in this paper are finite, undirected, and simple.
For notation not defined here, we refer the reader to
Section~\ref{sec:definition} and to \cite{bcn89,dkt}.

Let $\Gamma$ be a distance-regular graph with valency $k$, diameter $D \geq 3$,
and distinct eigenvalues $k = \theta_0 > \theta_1 > \cdots > \theta_D$. 
Let $b = \frac{b_1}{\theta_1 + 1}$. 
Then it is known that the smallest eigenvalue of each local graph of $\Gamma$ is at least $-b-1$ (see Lemma~\ref{localgraph}),  and that
$(\theta_1 + 1)(\theta_D + 1) < -b_1$
(see \cite[Theorem~3.6]{KPY2011}), which implies that $-\theta_D > b + 1$.

We are interested in whether there exists an upper bound for $-\theta_D$ that depends only on $b$ and $D$. For distance-regular graphs with classical parameters $(D,b,\alpha,\beta)$ and $b>0$, it is known that
$-\theta_D = 1 + b + b^2 + \cdots + b^{D-1}$. This leads us to the following problem.

\begin{problem}\label{pro:thetaD-bound}
Let $\Gamma$ be a distance-regular graph with valency $k$, diameter $D \geq 3$, and distinct eigenvalues
$k = \theta_0 > \theta_1 > \cdots > \theta_D$.
Let $b = \frac{b_1}{\theta_1 + 1}$.
Does there exist a function $f(D,b)$ such that $-\theta_D \leq f(D,b)$?
\end{problem}

Let $\Gamma$ be a distance-regular graph with valency $k$, smallest eigenvalue
$\theta_D$, and diameter $D \geq 3$. It is known that for a clique $C$ of
order $|C|$ in $\Gamma$, we have $|C| \leq 1 + \frac{k}{-\theta_D}$, and if equality holds, then $C$ is called a \emph{Delsarte clique}.
Let $C$ be a Delsarte clique in $\Gamma$. Then there exist numbers
$\phi_0, \phi_1, \ldots, \phi_{D-1}$ such that, for a vertex $x$ at distance $i$
from $C$, the number of vertices in $C$ at distance $i$ from $x$ is equal to
$\phi_i$, for $0 \leq i \leq D-1$.

A distance-regular graph $\Gamma$ is called \emph{geometric} if there exists a
set of Delsarte cliques $\mathcal{C}$ in $\Gamma$ such that each edge lies in a
unique clique $C \in \mathcal{C}$.
For a geometric distance-regular graph $\Gamma$, define $\tau_i$ as the number
of Delsarte cliques containing a vertex $y$ at distance $i-1$ from $x$, where
the vertices $x$ and $y$ are at distance $i$, for $1 \leq i \leq D$. The numbers
$\tau_i$ depend only on the distance $d(x,y) = i$, and not on the particular
choice of the vertices $x$ and $y$.

We now discuss bounds on the numbers $\phi_i$ defined above. The following result provides a bound on $\phi_1$ in terms of $b$ when $\frac{k}{-\theta_D}$ is sufficiently large. It is a direct consequence of Lemmas~\ref{2phi1} and~\ref{phi1bound}.

\begin{proposition}\label{phi1bound-geo}
Let $\Gamma$ be a geometric distance-regular graph with diameter $D \geq 3$,
distinct eigenvalues $k = \theta_0 > \theta_1 > \cdots > \theta_D$. Let $b = \frac{b_1}{\theta_1 + 1}$. If $\frac{k}{-\theta_D} \geq b^4 + 2b^3 + 3b^2 + b + 2$, then $\phi_1 \leq b^2 + b + 1$.
\end{proposition}

For general $i$, we conjecture the following.

\begin{conjecture}\label{conj:phi-bound}
Let $\Gamma$ be a geometric distance-regular graph with diameter $D$ and
distinct eigenvalues $k = \theta_0 > \theta_1 > \cdots > \theta_D$.
Let $b = \frac{b_1}{\theta_1 + 1}$. Then there exists a function $f(i,b)$ such
that, if $D \geq 2i + 1$, then $\phi_i \leq f(i,b)$.
\end{conjecture}

We will discuss constraints on $\phi_i$ in Section~\ref{sec:phibound}.

Concerning the parameters $\tau_i$, we also believe that they should be bounded by a function depending only on $i$ and $b$. We therefore propose the following conjecture.

\begin{conjecture}\label{conj:tau-bound}
Let $\Gamma$ be a geometric distance-regular graph with diameter $D \geq 3$ and distinct eigenvalues
$k = \theta_0 > \theta_1 > \cdots > \theta_D$.
Let $b = \frac{b_1}{\theta_1 + 1}$.
Then there exists a function $f(i,b)$ such that $\tau_i \leq f(i,b)$.
\end{conjecture}

Note that for geometric distance-regular graphs this conjecture generalizes Problem~\ref{pro:thetaD-bound}, which asks whether the smallest eigenvalue $\theta_D = -\tau_D$ is bounded by a function of $D$ and $b$.

Moreover, if Conjecture~\ref{conj:tau-bound} holds for $i = 2$, then for geometric distance-regular graphs, 
Lemma~\ref{tau2geqphi1} together with the relation $c_2 = \phi_1 \tau_2$ implies that $c_2$ is bounded by a function of $b$, as asked in the following Problem~\ref{c2bound}.



\begin{problem}[cf.~{\cite[Problem~66]{dkt}}]\label{c2bound}
For a coconnected distance-regular graph $\Gamma$, show that the intersection number $c_2$ is bounded above by a function of $\frac{b_1}{\theta_1 + 1}$. 
\end{problem}

Some progress on Problem~\ref{c2bound} has been made in
\cite{KLPY2025,TKCP2022}.

We now turn to our main result on geometric distance-regular graphs.

\begin{theorem}\label{main2}
Let $\Gamma$ be a geometric distance-regular graph with valency $k$, diameter $D \geq 3$, distinct eigenvalues
$k = \theta_0 > \theta_1 > \cdots > \theta_D$, and $c_2 \geq 2$.
Let $r = -\theta_D$, $\beta = \frac{k}{-\theta_D}$, and $b = \frac{b_1}{\theta_1 + 1}$.
Then one of the following holds: 
\begin{enumerate}
    \item $\phi_1 = 1$, and $\Gamma$ is locally the disjoint union of cliques of order $\beta$.
    \item $\phi_1 = \tau_2 = 2$, and $\Gamma$ is a Johnson graph. In particular, $b = 1$.
    \item $\phi_1 = \tau_2 \geq 3$, and $\Gamma$ is a Grassmann graph defined over the field $\mathbb{F}_{\phi_1-1}$. In particular, $b = \phi_1 - 1$.
    \item $\phi_1 = \tau_2 - 1 \geq 2$, and for each vertex $x$ of $\Gamma$, the local graph $\Delta_x$ at vertex $x$ is the $(\phi_1 - 1)$-clique extension of a $\frac{\beta}{\phi_1 - 1} \times r$-grid. In particular, if $r(\phi_1 - 1) > 2b^2 + 2b + 1$, then $\phi_1 \leq b + 1$.

    \item $ 2 \leq \phi_1 \leq \tau_2-1$, and $\beta < (r - \tau_2 + 1)(\phi_2 - \phi_1) + \phi_1$.
\end{enumerate}
\end{theorem}

Koolen and Bang~\cite{KB2010} proposed the following two conjectures concerning geometric distance-regular graphs.

\begin{conjecture}[cf.~{\cite[Conjecture 7.3]{KB2010}}]\label{conj7.3}
For a fixed integer $m \geq 2$, there are only finitely many coconnected geometric distance-regular graphs with smallest eigenvalue $-m$ and $\phi_1 \leq \tau_2 - 2$.
\end{conjecture}

\begin{conjecture}[cf.~{\cite[Conjecture 7.4]{KB2010}}]\label{conj7.4}
For a fixed integer $m \geq 2$, any geometric distance-regular graph with smallest eigenvalue $-m$, diameter $D \geq 3$ and $c_2 \geq 2$ is either a Johnson graph, a Grassmann graph, a Hamming graph, a bilinear forms graph, or the number of vertices is bounded above by a function of $m$.
\end{conjecture}

Note that Theorem~\ref{main2}, together with the constraints obtained for the parameters $\phi_i$ in Section~\ref{sec:phibound}, provides partial progress towards these two conjectures.

We conclude this section with an outline of the paper.
In the next section, we present preliminaries and basic definitions.
In Section~3, we establish restrictions on the parameters $\phi_i$.
In Section~4, we discuss the ELS property, which guarantees the existence of strongly regular subgraphs with desirable properties.
In Section~5, we derive explicit bounds on some geometric parameters.
In Section~6, we give a characterization of geometric distance-regular graphs with diameter $D\geq 3$ and $\tau_2=\phi_1\geq 2$.
In Section~7, we determine the local graphs of geometric distance-regular graphs that satisfy the dual Pasch axiom.
Finally, in Section~8, we present the proof of Theorem~\ref{main2}.

\section{Definitions and preliminaries}\label{sec:definition}
Let $\Gamma$ be a connected graph with vertex set $V(\Gamma)$. The \emph{distance} $d(x, y)$ between two vertices $x, y \in V(\Gamma)$ is the length of the shortest path in $\Gamma$ connecting $x$ and $y$. The \emph{diameter} $D$ of $\Gamma$ is the maximum distance between any two vertices of $\Gamma$. For two vertices $x, y \in V(\Gamma)$, we use $x \sim y$ to denote that $x$ is adjacent to $y$, and $x \not\sim y$ to denote that $x$ is not adjacent to $y$.
For each vertex $x \in V(\Gamma)$, let $\Gamma_i(x)$ denote the set of vertices in $\Gamma$ that are at distance $i$ from $x$, where $0 \leq i \leq D$. In addition, define $\Gamma_{-1}(x)=\Gamma_{D+1}(x)=\emptyset$. For the sake of simplicity, we denote $\Gamma_1(x)$ by $\Gamma(x)$.   For each vertex $x$ of $\Gamma$, the subgraph induced on $\Gamma(x)$ is called the \emph{local graph} of $\Gamma$ at vertex $x$, and we denote it by $\Delta_x$.  For a vertex $x$ of $\Gamma$, the  cardinality $|\Gamma(x)|$ of $\Gamma(x)$ is called the \emph{valency} of $x$ in $\Gamma$. In particular, $\Gamma$ is \emph{regular} with  valency $k$ if $k=|\Gamma(x)|$ holds for all $x\in V(\Gamma)$.  
For vertices $x$ and $y$ of $\Gamma$ at distance $i$ ($i = 1, 2, \dots, D$), we define the set $C_i(x, y)$ as $C_i(x, y) := \Gamma_{i-1}(x) \cap \Gamma(y)$. For vertices $u$ and $v$ of $\Gamma$ at distance $j$ ($j = 0, 1, \dots, D-1$), we define the set $B_j(u, v)$ as $B_j(u, v) := \Gamma_{j+1}(u) \cap \Gamma(v)$. For a vertex $x \in V(\Gamma)$ and a vertex subset $T \subseteq V(\Gamma)$, we define the distance $d(x, T)$ by $d(x, T) := \min\{d(x, y) \mid y \in T\}$.

A \emph{geodesic} between two vertices $x$ and $y$ in a graph $\Gamma$ is the shortest path connecting $x$ and $y$.
A set of vertices $T \subseteq V(\Gamma)$ is \emph{geodetically closed} if, for every pair of vertices $x, y \in T$, every shortest path between $x$ and $y$ lies entirely within $T$.

The \emph{adjacency matrix} $A=A(\Gamma)$ of $\Gamma$ is the matrix whose rows and columns are indexed by vertices of $\Gamma$ and the ($x, y$)-entry is $1$ whenever $x$ and $y$ are adjacent and $0$ otherwise. The \emph{eigenvalues} of $\Gamma$ are the eigenvalues of its adjacency matrix $A$.

A \emph{clique} of a graph $\Gamma$ is a set of mutually adjacent vertices of $\Gamma$. We sometimes also refer to a complete subgraph of $\Gamma$ as a clique.

The \emph{$m\times n$-grid} is the Cartesian product $K_m \Box K_n$ where $K_n$ is the complete graph on $n$ vertices. 
For a positive integer $s$ and a graph $\Gamma$, the \emph{$s$-clique extension of $\Gamma$}, denoted by $\tilde{\Gamma}$, is the graph obtained from $\Gamma$ by replacing each vertex $x$ with a complete graph $\tilde{X}$ having $s$ vertices, such that $u \in \tilde{X}$ is adjacent to $v\in \tilde{Y}$ $(\tilde{X}\neq \tilde{Y})$ in $\tilde{\Gamma}$ if and only if $x \sim y$ in $\Gamma$.

Let $\Pi=\{P_1,P_2,\dots,P_t\}$ be a partition of the vertex set of a graph $\Gamma$. Let $f_{ij}$ be the average number of neighbors in $P_j$ of a vertex in $P_i$, for $i,j=1,2,\dots,t$. The matrix $F=(f_{ij})$ is called the {\em quotient matrix} of $\Pi$. The partition $\Pi$ is called {\em equitable} if every vertex in $P_i$ has exactly $f_{ij}$ neighbors in $P_j$. For the definitions of the quotient matrix and equitable partition of a general symmetric real matrix, see \cite[Section~2.3]{BH2011}.

\subsection{Distance-regular graphs}  
A graph $\Gamma$ is called \emph{distance-regular} if there exist integers $b_i$ and $c_i$ ($0 \leq i \leq D$) such that, for each pair of vertices $x, y \in V(\Gamma)$ with $d(x, y) = i$, there are exactly  $c_i$ neighbors of $y$ in $\Gamma_{i-1}(x)$ and $b_i$ neighbors of $y$ in $\Gamma_{i+1}(x)$. Here, we define $b_D = c_0 = 0$.
Note that, in this case, $C_i(x, y)$ contains exactly $c_i$ vertices for $i = 1, 2, \ldots, D$, and $B_i(x, y)$ contains exactly $b_i$ vertices for $i = 0, 1, \ldots, D-1$.
In particular, every distance-regular graph is regular with valency $k = b_0$. We define $a_i := k - b_i - c_i$ for notational convenience.  
Note that $a_i = |\Gamma(y) \cap \Gamma_i(x)|$ holds for all pairs of vertices $x, y$ with $d(x, y) = i$ $(0 \leq i \leq D)$, and that the numbers $a_i$, $b_i$, and $c_i$ $(0 \leq i \leq D)$ are called the \emph{intersection numbers} of $\Gamma$. 

A graph $\Gamma$ is called a \emph{strongly regular graph} with parameters $(v, k, \lambda, \mu)$ if $\Gamma$ has $v\geq 2$ vertices, is $k$-regular and 
every pair of distinct vertices has 
$\lambda$ (resp. $\mu$) common neighbours depending on whether they are adjacent or not. Strongly regular graphs are precisely distance-regular graphs with diameter $D = 2$.

For each vertex $x \in V(\Gamma)$ and an integer $i$ ($0 \leq i \leq D$), there is a constant number of vertices at distance $i$ from $x$, denoted by $k_i$.
 That is, $k_i = |\Gamma_i(x)|$ for all $x \in V(\Gamma)$. Indeed, this follows by induction and by counting the number of edges between $\Gamma_i(x)$ and $\Gamma_{i+1}(x)$ in two different ways. In particular, we have $k_0 = 1$ and $k_{i+1} = \frac{b_i k_i}{c_{i+1}}$ for all $i = 0, 1, \ldots, D-1$.

A distance-regular graph $\Gamma$ with valency $k$ and diameter $D$ has exactly $D+1$ distinct eigenvalues (see \cite[p.~128]{bcn89}), which we denote by $k = \theta_0 > \theta_1 > \cdots > \theta_D$.
For an eigenvalue $\theta$ of $\Gamma$, the sequence $(u_i)_{i=0,1,...,D}$ = $(u_i(\theta))_{i=0,1,...,D}$
satisfying $u_0$ = $u_0(\theta)$ = $1$, $u_1$ = $u_1(\theta)$ = $\theta/k$, and
\begin{align}\label{standardseq}
c_i u_{i-1} + a_i u_i + b_i u_{i+1} = \theta u_i ~(i=1,2,\ldots,D-1)
\end{align}
is called the \emph{standard sequence} corresponding to the eigenvalue $\theta$ (see \cite[p.128]{bcn89}). For the convenience of the reader we give the following 
well-known lemma.

\begin{lemma}[cf.~{\cite[Proposition~4.4.1]{bcn89}}]\label{standard}
Let $\Gamma$ be a distance-regular graph, $\theta \neq k$ an eigenvalue of  $\Gamma$, say with multiplicity $m$ and let $(u_i)_i$ be the standard sequence corresponding to 
$\theta$. 
Then there exists a map $f: V(\Gamma) \rightarrow \mathbb{R}^m$, $x \mapsto \overline{x}$ such that the inner product
$\langle \overline{x}, \overline{y} \rangle = u_{d(x, y)}$. 
\end{lemma}

The map $f$ is called the \emph{standard representation} of $\Gamma$ with respect to $\theta$.

The following lemma is a well-known inequality for the smallest 
eigenvalues of local graphs of distance-regular graphs.

\begin{lemma}[cf.~{\cite[Theorem~4.4.3]{bcn89}}]\label{localgraph}
Let $\Gamma$ be a distance-regular graph with valency $k$, diameter $D \geq 3$, and distinct eigenvalues $k = \theta_0 > \theta_1 > \cdots > \theta_D$.
Let $b = \frac{b_1}{\theta_1 + 1}$.
Then, for a vertex $x$ of $\Gamma$, the smallest eigenvalue of the local graph $\Delta_x$ satisfies
\[
\theta_{\min}(\Delta_x) \geq -b - 1.
\]
\end{lemma}

\subsection{Distance-regular graphs with classical parameters}

A distance-regular graph $\Gamma$ is said to have \emph{classical parameters} $(D,b,\al,\beta)$
if the diameter of $\Gamma$ is $D$ and the intersection numbers of $\Gamma$ can be expressed as follows:
\begin{equation}\label{bi}
 b_i=([D]-[i])
(\beta-\al [i]),~~  0\le i\le D-1,
\end{equation}
\begin{equation}\label{ci}
c_i=[i](1+\al [i-1]),~~1\le i\le D,
\end{equation}
where \begin{equation}\label{coeff}
[j]= \begin{cases}
    \frac{b^j-1}{b-1} & \text{if } b\neq1, \\
    \binom{j}{1} & \text{if } b=1. \\
  \end{cases}\nonumber
\end{equation}
From \cite[Corollary 8.4.2]{bcn89}, we know that the eigenvalues of $\Gamma$ are 
\begin{equation}\label{eigen}
[D-i](\beta- \alpha[i])-[i]=\frac{b_i}{b^i}-[i], ~~0\leq i\leq D.
\end{equation}
We note that $c_2=(b+1)(\al+1)$ and that if $b\geq1$, then the eigenvalues $\theta_i=\frac{b_i}{b^i}-[i] (0\leq i\leq D)$ of $\Gamma$ are in the natural ordering, i.e., $k=\theta_0>\theta_1>\cdots>\theta_D$.

From the following result we know that the parameter $b$ of a distance-regular graph with classical parameters $(D,b,\alpha,\beta)$ is an integer.
\begin{proposition}[cf.~{\cite[Proposition~6.2.1]{bcn89}}]\label{b_integer}
Let $\Gamma$ be a distance-regular graph with classical parameters $(D,b,\alpha,\beta)$ and diameter $D\geq3$. Then $b$ is an integer such that $b\neq0,-1$. 
\end{proposition}

\subsection{Partial linear spaces}
An \emph{incidence structure} is a triple $(\mathcal{P}, \mathcal{L}, \mathcal{I})$ where $\mathcal{P}$ and $\mathcal{L}$ are nonempty disjoint sets and $\mathcal{I} \subseteq \mathcal{P} \times \mathcal{L}$. The elements of $\mathcal{P}$ and $\mathcal{L}$ are called \emph{points} and \emph{lines}, respectively. If $(p,\ell) \in \mathcal{I}$, then $p$ is said to be \emph{incident} with $\ell$. In this case, we also say that $\ell$ contains $p$ or that $p$ lies on $\ell$. An incidence structure $(\mathcal{P}, \mathcal{L}, \mathcal{I})$ is said to be \emph{finite} if both $\mathcal{P}$ and $\mathcal{L}$ are finite sets.
Two incidence structures $(\mathcal{P}, \mathcal{L}, \mathcal{I})$ and $(\mathcal{P}', \mathcal{L}', \mathcal{I}')$ are said to be  \emph{isomorphic} if  there exist bijections $\sigma: \mathcal{P} \to \mathcal{P}'$ and $\tau: \mathcal{L} \to \mathcal{L}'$ such that $(p, \ell) \in \mathcal{I}$ if and only if $(\sigma(p), \tau(\ell)) \in \mathcal{I}'$.

Two distinct points $p_1$ and $p_2$ are called \emph{collinear} if there exists a line $\ell$ containing both $p_1$ and $p_2$; otherwise, they are said to be \emph{noncollinear}. Two distinct lines $\ell_1$ and $\ell_2$ are said to \emph{intersect} if  they have a common incident point.

A \emph{partial linear space} is an incidence structure  such that each pair of distinct points is incident with at most one line. 
A \emph{linear space} is an incidence structure  such that each pair of distinct points is incident with exactly one line.

The \emph{point graph} $\Gamma$ of an incidence structure $(\mathcal{P}, \mathcal{L}, \mathcal{I})$ is the graph with vertex set $\mathcal{P}$ such that two distinct points are adjacent if and only if they lie on a common line. An incidence structure $(\mathcal{P}, \mathcal{L}, \mathcal{I})$ is said to be \emph{connected} if its point graph is connected.

If the incidence structure is connected, then for two distinct points $x, y$ we define their distance $d(x, y)$ to be their distance in the point graph $\Gamma$. For a point $x$ and a line $\ell$, we define the distance $d(x, \ell)$ by $d(x, \ell) := \min\{d(x, y)\mid y \text{ lies on } \ell\}$. For two points $x, y$ at distance $2$, we denote by $[x, y]$ the set of lines $\ell$ containing $x$ such that $d(\ell, y)=1$.

\subsection{Geometric distance-regular graphs}
The following lemma is called the Delsarte bound. For a proof, see \cite[Lemma 3]{LK2024+}.

\begin{lemma}[cf.~{\cite[Proposition~4.4.6]{bcn89}}]\label{Del}
Let $\Gamma$ be a distance-regular graph with diameter $D \geq 2$, valency $k$ and smallest eigenvalue $\theta_{\min}$, say with multiplicity $m$.
Let $(u_i)_i$ be the standard sequence corresponding to $\theta_{\min}$. 
Then the order $c$ of a clique $C$ in $\Gamma$ is bounded by $c \leq 1 + \frac{k}{-\theta_{\min}}$. 

Assume that $c = 1 + \frac{k}{-\theta_{\min}}$. 
Then each vertex of $\Gamma$ has distance at most $D-1$ to $C$. Let $x$ be a vertex with $d(x, C) = j$ for some $1\leq j \leq D-1$. Then the number of vertices $\phi_j(x)$ in $C$ at distance $j$ from $x$ satisfies:
\begin{equation}\label{decl}
(u_{j+1} - u_{j})\phi_j(x) = \left(1+ \frac{k}{-\theta_{\min}}\right) u_{j+1},
\end{equation}
and hence $\phi_j(x)$ does not depend on $x$, only on $j$.
\end{lemma}

A clique of $\Gamma$ with $1+\frac{k}{-\theta_{\min}}$ vertices is called a \emph{Delsarte clique}. 

A distance-regular graph $\Gamma$ is called \emph{geometric} if there exists a set of Delsarte cliques $\mathcal C$ in $\Gamma$ such that each edge lies in a unique $C \in \mathcal C$. In this case, we say that $\Gamma$ is a geometric distance-regular graph with respect to $\mathcal C$. Moreover, a partial linear space $X = (V(\Gamma), \mathcal C, \in)$ is naturally induced, whose point graph is $\Gamma$, and hence the elements of $\mathcal C$ are also called lines.

Let $\Gamma$ be a geometric distance-regular graph with respect to $\mathcal C$. Let $C$ be a Delsarte clique and $x$ be a vertex at distance $j$ from $C$ for $0 \leq j \leq D-1$. 
Define $\phi_j$ to be the number of vertices in $C$ at distance $j$ from $x$. 
The number $\phi_j$ $(0 \leq j \leq D-1)$ does not depend on the pair $x, C$, but only on their distance.
For two vertices $x$ and $y$ at distance $j$ for $1 \leq j \leq D,$ let $\tau_j$ be the number of cliques in $\mathcal C$ containing $x$ at distance $j-1$ from $y$. 
The number $\tau_j$ $(1 \leq j \leq D)$ does not depend on the pair $x, y$, but only on their distance.

Let $\Gamma$ be a geometric distance-regular graph. The numbers $\phi_i \ (i=0, 1, \ldots, D-1)$ and $\tau_j \ (j =1,2, \ldots, D)$ are called the \emph{geometric parameters} of $\Gamma$, and the array $\{\phi_0, \phi_1, \dots, \phi_{D-1}; \tau_1, \tau_2, \dots, \tau_D\}$ is called the \emph{geometric intersection array} of $\Gamma$. Note that $\phi_0 = 1$ and $\tau_1 = 1$ always hold.

We can express the intersection numbers of $\Gamma$ in terms of the geometric parameters as follows (see \cite[Lemma 4.1]{KB2010}):
\begin{equation}\label{eqgeom}
c_i = \tau_i \phi_{i-1} \ \ (i=1,2, \ldots, D), \ \ b_i= - (\theta_{\min} + \tau_i)\left(1+ \frac{k}{-\theta_{\min}} - \phi_i\right) \ \ (i = 1, 2, \ldots, D-1).
\end{equation}

For geometric distance-regular graphs with classical parameters, the following
result gives explicit expressions for $\phi_i$ and $\tau_j$.

\begin{lemma}[cf.~{\cite[Lemma~11]{LK2024+}}]\label{phi-CP}
Let $\Gamma$ be a distance-regular graph with classical parameters $(D, b, \alpha, \beta)$ such that $b \geq 2$ and $D \geq 3$, and assume that 
$\Gamma$ is geometric. Then the numbers $\phi_j$'s and $\tau_j$'s of Equation~(\ref{eqgeom}) satisfy
\begin{equation*} \phi_j = 1 + \alpha[j], \ \ (j =0, 1, \ldots, D-1), \ \ \tau_j = [j] \ \ (j =1,2, \ldots, D). \end{equation*}
\end{lemma}

In the following, we introduce four important classes of geometric distance-regular graphs with classical parameters: the Johnson graph $J(n, D)$, the Grassmann graph $J_q(n, D)$, the Hamming graph $H(D,e)$, and the bilinear forms graph $Bil(D \times e, q)$. Moreover, we characterize the Johnson graph and the Grassmann graph in Section~\ref{phi=tau}. For further details on these graphs, see \cite[Chapter~9]{bcn89} and \cite{dkt}.

The \emph{Johnson graph} $J(n, D)$ has as vertices the subsets of size $D$ of a set of size $n$. Two subsets are adjacent if and only if they differ in precisely one element. Note that $J(n, D)$ is isomorphic to $J(n, n - D)$. Therefore, in the following, we restrict to $n \geq 2D$.

The \emph{Grassmann graph} $J_q(n, D)$ has as vertices the $D$-dimensional subspaces of a vector space of dimension $n$ over the field $\mathbb{F}_{q}$. Two subspaces are adjacent if and only if they intersect in a $(D - 1)$-dimensional subspace. Note that $J_q(n, D)$ is isomorphic to $J_q(n, n - D)$. Again, we therefore restrict to $n \geq 2D$.

The \emph{Hamming graph} $H(D,e)$ has vertex set $X^D$, consisting of all words of length $D$ over an alphabet $X$ of size $e$. Two vertices are adjacent if and only if they differ in precisely one position. Hamming graphs are geometric distance-regular graphs with $\phi_1 = \tau_2 - 1 = 1$, and they satisfy item~(1) of Theorem~\ref{main2}.

The \emph{bilinear forms graph} $Bil(D \times e, q)$ has as its vertices all $D \times e$ matrices over the field $\mathbb{F}_{q}$, where two matrices are adjacent if and only if their difference has rank $1$. We assume $D \leq e$ throughout, so that $D$ is the diameter. The bilinear forms graph also admits an alternative description: its vertices are the $D$-dimensional subspaces of a $(D+e)$-dimensional vector space over $\mathbb{F}_{q}$ that intersect a fixed $e$-dimensional subspace trivially, and two vertices are adjacent whenever their intersection has dimension $D-1$. This description shows that the graph is isomorphic to a subgraph of the Grassmann graph $J_q(D+e, D)$. Bilinear forms graphs are geometric distance-regular graphs with $\phi_1 = \tau_2 - 1 = q$, and they satisfy item~(4) of Theorem~\ref{main2}.

The following lemmas give some key properties of geometric parameters, which play a crucial role in the characterization of geometric distance-regular graphs.

\begin{lemma}[cf.~{\cite[Lemma 4.2]{KB2010}}]\label{tau2geqphi1}
Let $\Gamma$ be a geometric distance-regular graph with diameter $D \geq 2$ and intersection number $c_2 \geq 2$. Then the following holds:  
\begin{itemize}  
    \item[(i)] $\tau_2 \geq \phi_1$.  
    \item[(ii)] $\Gamma$ contains an induced quadrangle.  
\end{itemize}
\end{lemma}

\begin{lemma}[cf.~{\cite[Proposition 4.3]{KB2010}}]\label{Dbound}
Fix an integer $m \geq 2$. Let $\Gamma$ be a geometric distance-regular graph with smallest eigenvalue $-m$, diameter $D \geq 2$. If $c_2 \geq 2$, then $D < m^2$.
\end{lemma}

\begin{lemma}[cf.~{\cite[Theorem 5.5]{BHK2007}}]\label{phiincrease}
Let $\Gamma$ be a geometric distance-regular graph with valency $k$ and diameter $D \geq 2$. If $\phi_1 > 1$, then
\[
1 < \phi_1 < \cdots < \phi_{D-1}.
\]  
In particular, $D \leq s$, where $s + 1$ is the order of a Delsarte clique.

\end{lemma}

\begin{lemma}[cf.~{\cite[Lemma 4.1]{bang2018}}]\label{tauincrease}
Fix an integer $m \geq 2$. Let $\Gamma$ be a geometric distance-regular graph with diameter $D \geq 2$ and smallest eigenvalue $-m$. If $\phi_1 \geq 2$, then
\[
2 \leq \phi_1 \leq \tau_2 < ... < \tau_D = m.
\]
\end{lemma}

\begin{lemma}\label{tauinequality}
    Let $\Gamma$ be a geometric distance-regular graph with diameter $D \geq 2$.  
    Let $i$ and $j$ be positive integers such that $i + j \leq D$. Then 
    \[
    \tau_i + \tau_j \leq \tau_D.
    \]
\end{lemma}
{\bf Proof.}
Let $\Gamma$ be a geometric distance-regular graph with respect to $\mathcal C$.  
Let $x$ and $y$ be a pair of vertices at distance $i + j$, and let $z$ be a vertex such that $d(x, z) = i$ and $d(y, z) = j$.  
Let $\mathcal{A} = \{C_1, C_2, \dots, C_{\tau_i}\}$ be the set of Delsarte cliques in $\mathcal{C}$ containing $z$ at distance $i - 1$ from $x$, and let $\mathcal{B} = \{C'_1, C'_2, \dots, C'_{\tau_j}\}$ be the set of Delsarte cliques in $\mathcal{C}$ containing $z$ at distance $j - 1$ from $y$.  

For each clique $C \in \mathcal{A}$, there exists a vertex $u \in C$ at distance $i - 1$ from $x$.  
Since $d(x, y) = i + j$ and $d(y, z) = j$, the vertex $u$ is at distance $j + 1$ from $y$.  
Therefore, $u$ cannot lie on any clique $C' \in \mathcal{B}$.  
Hence, $\mathcal{A} \cap \mathcal{B} = \emptyset$. Since $z$ lies on exactly $\tau_D$ Delsarte cliques in $\mathcal{C}$, we conclude that $\tau_i + \tau_j \leq \tau_D$.
\qed

\begin{lemma}[cf.~{\cite[Lemma 5.2]{BHK2007}}]\label{2phi1}
Let $\Gamma$ be a distance-regular graph with valency $k$ and diameter $D \geq 2$, containing a Delsarte clique $C$ of order $s + 1$. Then the following hold:
\begin{enumerate}
    \item Let $i$ and $j$ be positive integers such that $i + j \leq D - 1$. Then $\phi_i + \phi_j \leq s + 1$.
    \item If $\phi_1 > \frac{s + 1}{2}$, then $D = 2$.
    \item Suppose $\phi_1 \geq 2$. Then for all $x \in C$, the local graph $\Delta_x$ is connected, and its second largest eigenvalue equals $s - \phi_1$. Moreover, if $D = 2$, then $\Gamma$ has eigenvalues $k, s - \phi_1, -\frac{k}{s}$.
\end{enumerate}
\end{lemma}

\subsection{Neumaier's bound for parameters of partial geometries}\label{neumaierSRG}

A \emph{partial geometry} $pg(s, t, \alpha)$ is a partial linear space $X=(\mathcal{P}, \mathcal{L}, \mathcal{I})$ that satisfies the following conditions:
\begin{enumerate}
\item Every line contains exactly $s+1 \geq 2$ points,
\item Every point lies on exactly $t+1 \geq 2$ lines,

\item For a nonincident pair $(x, \ell)$ of a point $x$ and line $\ell$, there are exactly $\alpha \geq 1$ points on $\ell$ that are collinear with $x$.
\end{enumerate}
The \emph{dual} of a partial geometry is obtained by interchanging the roles of points and lines. Note that the dual of a $pg(s, t, \alpha)$ is a partial geometry $pg(t, s, \alpha)$.
The point graph of a partial geometry is a geometric strongly regular graph, and conversely, a geometric strongly regular graph can be the point graph of a partial geometry.

The following theorem, given by Neumaier \cite{neumaier1979}, imposes constraints on the parameters of partial geometries, and consequently, on the parameters of geometric 
strongly regular graphs.

\begin{theorem}[cf.~{\cite[Theorem~4.5]{neumaier1979}}]\label{SRGbound}
If a $pg(s,t,\alpha)$ satisfies $\alpha \leq t-1$, then
$$
s \leq (t - \alpha + 1)^2 (2\alpha - 1),
$$
and equality implies $\alpha = 1$ or $t = 2\alpha$.
\end{theorem}

The dualization process described earlier allows us to interchange $s$ and $t$ in this result, and in fact we have stated the dual of the result as it was originally stated in \cite{neumaier1979}.

\subsection{Designs}
Let $v$, $k$, $t$, and $\lambda$ be integers satisfying $v \geq k \geq t \geq 0$ and $\lambda \geq 1$. 
Let $\mathcal{P}$ be a set of $v$ elements and $\mathcal{B}$ a subset of $\binom{\mathcal{P}}{k}$.
The pair $X = (\mathcal{P}, \mathcal{B})$ is called a \emph{$t$-$(v, k, \lambda)$-design} if each $T \in \binom{\mathcal{P}}{t}$ lies in exactly $\lambda$ elements of $\mathcal{B}$. 
The sets $\mathcal{P}$ and $\mathcal{B}$ are called the \emph{point set} and the \emph{block set} of $X$, respectively. The cardinality of $\mathcal{P}$ is called the \emph{size} of the design $X = (\mathcal{P}, \mathcal{B})$.

Let $X = (\mathcal{P}, \mathcal{B})$ be a $t$-$(v, k, \lambda)$-design. 
A \emph{subdesign} is a pair $(\mathcal{P}', \mathcal{B}')$ such that $\mathcal{P}' \subseteq \mathcal{P}$, $\mathcal{B}' \subseteq \mathcal{B}$, and $(\mathcal{P}', \mathcal{B}')$ 
is a $t$-$(v', k, \lambda')$-design. Clearly, the parameters of $(\mathcal{P}', \mathcal{B}')$ satisfy $v' \leq v$ and $\lambda' \leq \lambda$.

The following lemma, due to Rosa, establishes a relationship between the sizes of a $2$-$(v, k, 1)$-design and its subdesigns. For more properties of $2$-$(v, k, 1)$-designs, 
see \cite{Stin}.

\begin{lemma}[cf.~{\cite[Theorem~9.43]{Stin}}]\label{subdesign}
Let $X = (\mathcal{P}, \mathcal{B})$ be a $2$-$(v, k, 1)$-design with a subdesign of size $v' <v$. Then $v \geq (k - 1)v' + 1$.
In particular, if $v > k$, then $v \geq (k-1)k + 1$, since every block defines a $2$-$(k, k, 1)$-design.
\end{lemma}

\section{Restrictions on $\phi_i$}\label{sec:phibound}
In this section, we establish restrictions on the parameters $\phi_i$. To this end, we first derive restrictions for the more general quantities $\psi_j(x,C)$.

Let $\Gamma$ be a distance-regular graph. Let $C$ be a clique, and let $x$ be a vertex at distance $j$ from $C$, where $1 \leq j \leq D$. We define $\psi_j(x,C)$ to be the number of vertices of $C$ that are at distance $j$ from $x$.

The following lemma provides an upper bound for $\psi_1(x,C)$ in terms of
$b = \frac{b_1}{\theta_1 + 1}$.

\begin{lemma}\label{phi1bound}
Let $\Gamma$ be a distance-regular graph with valency $k$, diameter $D \geq 3$, and second largest eigenvalue $\theta_1$. Let $b = \frac{b_1}{\theta_1 + 1}$. Suppose that $C$ is a clique of order $s+1$, and that $x$ is a vertex at distance $1$ from $C$. If $s \geq b^4 + 2b^3 + 3b^2 + b + 2$ and $\psi_1(x,C) < s+1$, then either $\psi_1(x,C) \leq b^2 + b + 1$ or $\psi_1(x,C) \geq s - b^2 + 1$.
\end{lemma}
{\bf Proof.}
For a vertex $y \in C \cap \Gamma(x)$, let $\Delta_y$ denote the local graph at $y$.
By Lemma~\ref{localgraph}, the smallest eigenvalue of $\Delta_y$ satisfies $\theta_{\min}(\Delta_y) \geq -b - 1.$

Now consider $\Delta_y$. The set $C \setminus \{y\}$ is a clique in $\Delta_y$ of
order $s \geq b^4 + 2b^3 + 3b^2 + b + 2$, and the vertex $x \in V(\Delta_y)$ has
$\psi_1(x,C) - 1$ neighbours in $C \setminus \{y\}$. Therefore, by
\cite[Proposition~1.2]{YK}, we obtain either
$\psi_1(x,C) \leq b^2 + b + 1$ or $\psi_1(x,C) \geq s - b^2 + 1$.
\qed

\vspace{0.2cm}

The following result, motivated by \cite[Lemma~3.1]{GKP2021} and
\cite[Proposition~1.2]{YK}, gives constraints on the size of the intersection
of two large cliques. It will be used to prove that, in the situation described
in Theorem~\ref{main2}\,(4), we have $\phi_1 \leq b + 1$.
\begin{lemma}\label{clique-int-bound}
Let $m \geq 1$ be an integer and $\Gamma$ a graph with smallest eigenvalue
$\theta_{\min} \geq -m$. Let $C_1$ and $C_2$ be two intersecting cliques in $\Gamma$
with $\min\{|C_1|, |C_2|\} = c$. If
\[
c > 2m^2-2m+1,
\]
then either $|C_1 \cap C_2| \leq m-1$, or $|C_1 \cap C_2| \geq c - m + 1$.
\end{lemma}
{\bf Proof.}
Without loss of generality, we may assume that $|C_1| = |C_2| = c$.
Let $\gamma = |C_1 \cap C_2|$. Assume that $m \leq \gamma \leq c-m$.
Let $H$ be the subgraph of $\Gamma$ induced on $C_1 \cup C_2$.  Observe that $\pi=\{\, C_1 \cap C_2,\ (C_1 \cup C_2)\setminus (C_1 \cap C_2)\,\}$
is an equitable partition of the vertex set of $H$, and the corresponding quotient
matrix $Q$ is given by
\[
Q=
\begin{pmatrix}
\gamma-1 & 2(c-\gamma)\\
\gamma & c-\gamma-1
\end{pmatrix}.
\]

By \cite[Lemma~2.3.1]{BH2011}, the smallest eigenvalue of $Q$ is at least $-m$.
Hence $\det(mI+Q)\geq 0$, which yields
\[
(m+\gamma-1)(m+c-\gamma-1)\geq 2\gamma(c-\gamma).
\]
Equivalently,
\[
(m-1)^2 \geq c(\gamma-m+1)-\gamma^2.
\]

Let $f(\gamma)=c(\gamma-m+1)-\gamma^2$. Then
\[
(m-1)^2 \geq f(\gamma)
\geq \min\{f(m),f(c-m)\}
= c-m^2.
\]
Therefore,
$c \leq 2m^2-2m+1$. Hence, the lemma holds.
\qed

\vspace{0.2cm}

We now consider $\psi_j(x,C)$ in general.
The following Lemma~\ref{psibound} provides bounds on $\psi_j(x,C)$.

\begin{lemma}\label{psibound}
Let $\Gamma$ be a distance-regular graph with valency $k$ and diameter $D \geq 2$, and let $\theta$ be an eigenvalue of $\Gamma$ with corresponding standard sequence $(u_i)_i$. Suppose that $C$ is a clique of order $s+1$, and that $x$ is a vertex at distance $j$ from $C$, where $1 \leq j \leq D-1$. If $\psi_j(x,C) < s+1$, then the following inequality holds:

\[
\begin{aligned}
\bigl(\psi_j(x,C) u_1 (u_j - u_{j+1}) 
    &- (1 - u_1) u_{j+1} \bigr)^2  \\
&\geq (su_1 + 1) \Bigl( \psi_j(x,C) u_1 (u_j - u_{j+1})^2 
    - (1 - u_1)\bigl(u_1 - u_{j+1}^2\bigr) \Bigr).
\end{aligned}
\]
\end{lemma}
{\bf Proof.}
Let $A = \{\, y \in C \mid d(x,y) = j \,\}$ and $B = \{\, y \in C \mid d(x,y) = j+1 \,\}$. As $\psi_j(x,C) < s+1$, it follows that $B \neq \emptyset$. For convenience, we write $\psi_j = \psi_j(x,C)$.

Let $\sigma : V(\Gamma) \to \mathbb{R}^m$ denote the standard representation of $\Gamma$ with respect to $\theta$, where $m$ is the multiplicity of $\theta$. Let $M$ be the Gram matrix of the set $\{\overline{y} = \sigma(y) \mid y \in \{x\} \cup C \}$; that is, the $(y_1,y_2)$-entry of $M$ is the inner product $\langle \overline{y_1}, \overline{y_2} \rangle = u_{d(y_1,y_2)}$. 
Observe that $\pi = \{\{x\}, A, B\}$ is an equitable partition of $M$, and the quotient matrix $Q$ of $M$ with respect to $\pi$ is given by:

\begin{equation*}
		Q=
		\begin{pmatrix}
			1 & \psi_j u_j & (s+1 - \psi_j) u_{j+1} \\
            u_j & 1 + (\psi_j - 1) u_1 & (s+1 - \psi_j) u_1 \\
            u_{j+1} & \psi_j u_1 & 1 + (s - \psi_j) u_1
		\end{pmatrix}.
		\end{equation*}

\allowdisplaybreaks 
  Since $M$ is a positive semi-definite matrix, it follows that $det(Q)\geqslant 0$. Computing the determinant of $Q$, we have
\begin{align*}
\det(Q) 
&= \det\begin{pmatrix}
1 & \psi_j u_j & (s+1 - \psi_j) u_{j+1} \\[4pt]
u_j & 1 + (\psi_j - 1) u_1 & (s+1 - \psi_j) u_1 \\[4pt]
u_{j+1} & \psi_j u_1 & 1 + (s - \psi_j) u_1
\end{pmatrix} \\
&= \det\begin{pmatrix}
1 & \psi_j u_j & (s+1 - \psi_j) u_{j+1} \\[4pt]
u_j - u_{j+1} & 1 - u_1 &  u_1- 1 \\[4pt]
u_{j+1} & \psi_j u_1 & 1 + (s - \psi_j) u_1
\end{pmatrix} \\
&= \det\begin{pmatrix}
1 & \psi_j u_j & (s+1) u_{j+1} +  \psi_j (u_j - u_{j+1}) \\[4pt]
u_j - u_{j+1} & 1 - u_1 &  0 \\[4pt]
u_{j+1} & \psi_j u_1 & 1 + s  u_1
\end{pmatrix} \\
&= \det\begin{pmatrix}
u_1 & \psi_j u_1 u_j & (s+1) u_1 u_{j+1} +  \psi_j u_1 (u_j - u_{j+1}) \\[4pt]
u_j - u_{j+1} & 1 - u_1 &  0 \\[4pt]
u_{j+1} & \psi_j u_1 & 1 + s  u_1
\end{pmatrix} \\
&= \det\begin{pmatrix}
u_1 - u_{j+1}^2 & \psi_j u_1(u_j - u_{j+1})  & \psi_j u_1(u_j - u_{j+1}) -   (1 - u_1)u_{j+1} \\[4pt]
u_j - u_{j+1} & 1 - u_1 &  0 \\[4pt]
u_{j+1} & \psi_j u_1 & 1 + s  u_1
\end{pmatrix}. \\
\end{align*}
Therefore, we obtain the following inequality:
\[
\bigl(\psi_j u_1 (u_j - u_{j+1}) - (1 - u_1) u_{j+1} \bigr)^2 
\;\geq\; 
(s u_1 + 1) \Bigl( \psi_j u_1 (u_j - u_{j+1})^2 - (1 - u_1)\bigl(u_1 - u_{j+1}^2 \bigr) \Bigr).
\]

This completes the proof of Lemma~\ref{psibound}. 
\qed

\vspace{0.2cm}

Here we remind the reader that although Lemma~\ref{psibound} holds for every eigenvalue of $\Gamma$, motivated by the insight obtained from Lemma~\ref{phi1bound}, we will mainly consider the standard sequence $(u_i)_i$ corresponding to the second largest eigenvalue. Indeed, we expect that the second largest eigenvalue is likely to provide the strongest bounds.

The following is an immediate corollary of Lemma~\ref{psibound}. 
Note that for a Delsarte clique $C$ of order $s+1$ and a vertex $x$ at distance
$j$ from $C$, where $1 \leq j \leq D-1$, it follows from Lemma~\ref{Del} that
$\psi_j(x,C)$ is independent of the choice of $x$. Hence we may denote
$\phi_j := \psi_j(x,C)$, which is well defined and satisfies $\phi_j < s+1$.

\begin{corollary}\label{psibound-DC}
Let $\Gamma$ be a distance-regular graph with valency $k$ and diameter $D \geq 2$, and let $\theta_1$ be the second largest eigenvalue of $\Gamma$ with corresponding standard sequence $(u_i)_i$. Suppose that $C$ is a Delsarte clique of order $s+1$. Then, for $1 \leq j \leq D-1$, the following inequality holds:
\[
\bigl(\phi_j u_1 (u_j - u_{j+1}) - (1 - u_1) u_{j+1} \bigr)^2 
\;\geq\; 
(s u_1 + 1) \Bigl( \phi_j u_1 (u_j - u_{j+1})^2 - (1 - u_1)\bigl(u_1 - u_{j+1}^2 \bigr) \Bigr).
\]
In particular, the above inequality holds for geometric distance-regular graphs.
\end{corollary}

\begin{remark}
Corollary~\ref{psibound-DC} implies that if
$\phi_j u_1 (u_j - u_{j+1})^2 \gg (1 - u_1)\bigl(u_1 - u_{j+1}^2\bigr)$,
then $\phi_j$ is close to $s$.
\end{remark}


We make the above remark explicit for geometric distance-regular graphs with
classical parameters in the following corollary.

\begin{corollary}\label{psibound-CP}
Let $\Gamma$ be a geometric distance-regular graph with classical parameters $(D, b, \alpha, \beta)$, valency $k$, and diameter $D \geq 3$, and let the distinct eigenvalues of $\Gamma$ be $k = \theta_0 > \theta_1 > \cdots > \theta_D$.  Let $(u_i)_i$ denote the standard sequence corresponding to $\theta_1$. Then the following inequalities hold:
\begin{enumerate}
\item $\left( \phi_1 - \frac{b u_2}{u_1} \right)^2 \geq \left( \beta + \frac{1}{u_1} \right) \left( \phi_1 - b^2 - 2b - 1 + \frac{1}{u_1} \right).$

\item $\left(\phi_2 - \frac{b^2 u_3}{u_1}\right)^2 \geq \left(\beta + \frac{1}{u_1}\right)\left(\phi_2 - b^4 - 2b^3 - b^2\left(1 + \frac{2u_2-1}{u_1}\right) - 1 + \frac{1}{u_1}\right).$
\end{enumerate}
In particular, if $b \geq 2$, then the following inequalities hold:
\begin{enumerate}[resume]
\item $\left( \phi_1 + \frac{7b^2}{2} \right)^2 \geq  \beta \left( \phi_1 - b^2 - 2b - 1 \right)$.

\item $\left(\phi_2 + \frac{7b^3}{2}\right)^2  \geq \beta \left( \phi_2 - b^4 - \frac{ 11b^3}{2} - b^2 - 1\right)$.
\end{enumerate}
\end{corollary}
{\bf Proof.}
Define $r = -\theta_D$.
By equations~\eqref{bi} and~\eqref{eigen}, we obtain 
$ r = -\theta_D = -[D] $, 
$ k = b_0 = r\beta $, 
$ b_1 = (r-1)(\beta - \alpha) $, 
and 
$ \theta_1 = \frac{b_1}{b} - 1 = \frac{(r-1)(\beta - \alpha) - b}{b} $.

By Lemma~\ref{Del}, each Delsarte clique has $\frac{k}{-\theta_D} + 1 = \beta + 1$ vertices. Applying Corollary~\ref{psibound-DC} with $j = 1, 2$, we obtain 
\begin{equation}\label{phi1-ineq}
\bigl(\phi_1 u_1 (u_1 - u_2) - (1 - u_1) u_2 \bigr)^2 
\;\geq\; 
(\beta u_1 + 1)\Bigl( \phi_1 u_1 (u_1 - u_2)^2 - (1 - u_1)\bigl(u_1 - u_2^2\bigr) \Bigr),
\end{equation}
and
\begin{equation}\label{phi2-ineq}
\bigl(\phi_2 u_1 (u_2 - u_3) - (1 - u_1) u_3 \bigr)^2 
\;\geq\; 
(\beta u_1 + 1) \Bigl( \phi_2 u_1 (u_2 - u_3)^2 - (1 - u_1)\bigl(u_1 - u_3^2 \bigr) \Bigr).
\end{equation}

By \cite[Lemma~4.6]{KAGS2024}, there exists a real number $c$ such that 
\begin{equation}\label{ui-eq}
u_i = \frac{1}{b} u_{i-1} + c
\qquad (1 \leq i \leq D).
\end{equation}

By equations~\eqref{standardseq} and~\eqref{ui-eq}, we obtain  
\[
\frac{1}{b} + c 
= u_1 
= \frac{\theta_1}{k} 
= \frac{(r-1)(\beta - \alpha) - b}{br\beta}, 
\]
and hence  
\[
c 
= -\frac{\beta + r\alpha + b - \alpha}{br\beta}.
\]

By equation~\eqref{ui-eq}, we have  
\begin{equation}\label{u1u2-diff}
u_1 - u_2 = \frac{1}{b}(1 - u_1).
\end{equation}
Hence,
\begin{align}
u_2^2 - u_1
&= \left( \frac{1 - u_1}{b} - u_1 \right)^2 - u_1   \nonumber \\
&= \frac{(1-u_1)^2}{b^2} - \frac{2u_1(1-u_1)}{b} - u_1(1-u_1).
\label{u2square}
\end{align}

Since $\theta_1 < k$, it follows that 
$1 - u_1 = 1 - \frac{\theta_1}{k} \neq 0$. Moreover, as $D \geq 3$, we have 
$\theta_1 > 0$ by \cite[Corollary~3.5.4]{bcn89}, and hence $u_1 > 0$. 
By Proposition~\ref{b_integer}, we also have $b \neq 0$.
Substituting equations~\eqref{u1u2-diff} and~\eqref{u2square} into~\eqref{phi1-ineq} 
and simplifying, we obtain
\[
\left( \phi_1 - \frac{b u_2}{u_1} \right)^2 \geq \left( \beta + \frac{1}{u_1} \right) \left( \phi_1 - b^2 - 2b - 1 + \frac{1}{u_1} \right),
\]
which shows Item~(1).

By equation~\eqref{ui-eq} and \eqref{u1u2-diff}, we have   
\begin{equation}\label{u2u3-diff}
u_2 - u_3 = \frac{u_1 - u_2}{b} = \frac{1 - u_1}{b^2}.
\end{equation}

By equation~\eqref{u2square} and \eqref{u2u3-diff}, we have 
\begin{align}
u_3^2 - u_1 
&= \left( \frac{1 - u_1}{b^2} - u_2 \right)^2 - u_1 \nonumber\\
&= \frac{(1 - u_1)^2}{b^4} - \frac{2u_2(1 - u_1)}{b^2} + u_2^2 - u_1 \nonumber\\
&= \frac{(1 - u_1)^2}{b^4} - \frac{2u_2(1 - u_1)}{b^2} + \frac{(1 - u_1)^2}{b^2} - \frac{2u_1(1 - u_1)}{b} - u_1(1 - u_1) \nonumber\\
&= \frac{(1 - u_1)^2}{b^4}
   + \frac{(1 - u_1 - 2u_2)(1 - u_1)}{b^2}
   - \frac{2u_1(1 - u_1)}{b}
   - u_1(1 - u_1).\label{u3square}
\end{align}

Substituting equations~\eqref{u2u3-diff} and~\eqref{u3square} into~\eqref{phi2-ineq} 
and simplifying, we obtain
\[
\left(\phi_2 - \frac{b^2 u_3}{u_1}\right)^2 \geq \left(\beta + \frac{1}{u_1}\right)\left(\phi_2 - b^4 - 2b^3 - b^2\left(1 + \frac{2u_2-1}{u_1}\right) - 1 + \frac{1}{u_1}\right),
\]
which shows Item~(2).

If $b \geq 2$, then by \cite[Lemma~8]{LK2024+} we have $\alpha \geq 0$. Moreover, by \cite[Proposition~1]{JV2017}, $\beta \geq \frac{\alpha (r-1)}{b} + 1 \geq 1.$ 
Therefore, we obtain  
\begin{align*}
-c 
&= \frac{\beta + r\alpha + b - \alpha}{br\beta} 
 = \frac{\beta + (r-1)\alpha + b}{br\beta} \\
&\leq \frac{\beta + b\beta + b}{br\beta}
 \leq \frac{1}{b}\left( \frac{b+1}{r} + \frac{b}{r} \right) \\
&\leq \frac{1}{b} \cdot \frac{2b+1}{b^{2}+b+1}
 \leq \frac{5}{7b}. 
\end{align*}
Hence, by equation~\eqref{ui-eq},  
\begin{equation}\label{u1-lbound}
u_1 = \frac{1}{b} + c \geq \frac{1}{b} - \frac{5}{7b} = \frac{2}{7b}.
\end{equation}

Now, if $u_2 \geq 0$, then $\left(\phi_1 - \frac{b u_2}{u_1}\right)^2 \leq \phi_1^2.$
If $u_2 < 0$, then by \cite[Theorem~4.4.1]{bcn89} we have $|u_2| \leq 1$, and thus  $\left(\phi_1 - \frac{b u_2}{u_1}\right)^2  
   \leq \left(\phi_1 + \frac{b}{u_1}\right)^2
   \leq \left(\phi_1 + \frac{7b^2}{2}\right)^2.$
   Therefore, in all cases, we have
\[
\left(\phi_1 - \frac{b u_2}{u_1}\right)^2  
   \leq \left(\phi_1 + \frac{7b^2}{2}\right)^2.
\]
Hence Item~(3) follows immediately from Item~(1).

If $u_3 \geq 0$, then $\left(\phi_2 - \frac{b^2 u_3}{u_1}\right)^2 \leq \phi_2^2$. If $u_3 < 0$, then by \cite[Theorem~4.4.1]{bcn89} we have $|u_3| \leq 1$, and thus $\left(\phi_2 - \frac{b^2 u_3}{u_1}\right)^2 \leq \left(\phi_2 + \frac{b^2}{u_1}\right)^2 \leq \left(\phi_2 + \frac{7b^3}{2}\right)^2$. Therefore, in all cases, we have
\[
\left(\phi_2 - \frac{b^2 u_3}{u_1}\right)^2 \leq \left(\phi_2 + \frac{7b^3}{2}\right)^2.
\]
Moreover,
\begin{align*}
&\quad \ \phi_2 - b^4 - 2b^3 - b^2\left(1 + \frac{2u_2-1}{u_1}\right) - 1 + \frac{1}{u_1} \\
&= \phi_2 - b^4 - 2b^3 - b^2 - 1 + \frac{1}{u_1}
      + b^2 \left(\frac{1-2u_2}{u_1}\right) \\
&\geq \phi_2 - b^4 - 2b^3 - b^2 - 1 - \frac{ b^2}{u_1}  \\
&\geq \phi_2 - b^4 - 2b^3 - b^2 - 1 - \frac{ 7b^3}{2} \\
&\geq \phi_2 - b^4 - \frac{ 11b^3}{2} - b^2 - 1.
\end{align*}
and therefore Item~(4) follows immediately from Item~(2).
\qed

\begin{remark}
\begin{enumerate}
    \item Corollary~\ref{psibound-CP}(1) and (3) show that if $\phi_1 \gg b^2$, then $\beta - \phi_1$ is small. However, since $D \geq 3$, Lemma~\ref{2phi1} implies that $\phi_1 \leq \frac{\beta+1}{2}$. Therefore, the quantity $|\phi_1 - b^2|$ should be small. This is also shown in Lemma~\ref{phi1bound}, with explicit bounds.

    \item Similarly, Corollary~\ref{psibound-CP}(2) and (4) show that if $\phi_2 \gg b^4$, then $\beta - \phi_2$ is small. However, when $D \geq 5$, Lemma~\ref{2phi1} implies that $\phi_2 \leq \frac{\beta+1}{2}$. Therefore, the quantity $|\phi_2 - b^4|$ should be small. The precise computation relies on the particular choice of the parameters $(D,b,\alpha,\beta)$.
    \item We believe that the above facts hold for general distance-regular graphs.
\end{enumerate}
\end{remark}

\section{The ELS property and strongly regular subgraphs}
Let $\Gamma$ be a geometric distance-regular graph with respect to a set $\mathcal C$ of Delsarte cliques. 
We say that $\Gamma$ has the \emph{equal line set} (ELS) property if the following condition holds:
In the partial linear space $X = (V(\Gamma), \mathcal C, \in)$, for each pair $x \in V(\Gamma)$ and $\ell \in \mathcal C$ such that $d(x,\ell)=1$, and for every pair of distinct vertices $y_1$ and $y_2$ on $\ell$ at distance $2$ from $x$, the sets of lines $[x,y_1]$ and $[x,y_2]$ are equal.

The following theorem provides a sufficient condition for a geometric distance-regular graph to satisfy the ELS property. In \cite[Theorem 22]{LK2024+}, we presented a similar result
for distance-regular graphs with classical parameters. We follow the proof of  \cite[Theorem 22]{LK2024+}, with some minor adjustments. The proof is based on the ideas of the proof of 
\cite[Lemma 2.10]{metsch}.

\begin{theorem}\label{geometric[]}
Let $\Gamma$ be a geometric distance-regular graph with respect to $\mathcal C$, having valency $k$, smallest eigenvalue $\theta_{\min}$, and diameter $D \geq 3$. Let $r = -\theta_{\min}$ and $\beta = \frac{k}{-\theta_{\min}}$.  
If $\beta \geq (r -\tau_2 + 1)(\phi_2-\phi_1)+\phi_1$, then $\Gamma$ satisfies the ELS property.

\end{theorem}
{\bf Proof.}
Let $\ell$ be a line at distance 1 from $x$.
Assume that there are vertices $y_1, y_2$ at distance 2 from $x$ that are both on $\ell$ such that $[x, y_1] \neq [x, y_2]$. Put, for $i=1,2$,  $S_i := \Gamma_i(x) \cap \ell$ and 
$n_i := |S_i|$.  By Lemma~\ref{Del}, we have $|\ell| = \beta +1$, $n_1 = \phi_1 $ and $n_2 = |\ell| - n_1 = \beta -\phi_1+1$, as $\Gamma$ is geometric.
This means that there are $\phi_1$  lines through $x$ that intersect $\ell$, as two distinct lines intersect in at most one point and $n_1 = \phi_1$. 
Those lines are in $T := [x, y_1] \cap [x, y_2]$. Let $t := |T|$. 

As $|[x, y_1]|= |[x, y_2]| = \tau_2$, and $[x, y_1] \neq [x, y_2]$, we have $\phi_1 = n_1 \leq t \leq \tau_2-1$. 
Because there are exactly $r$ lines through $x$, there are exactly $r -t$ lines through $x$ that are not in $T$. 
Let $u \in \Gamma(x)$ such that $d(y_1, u) =3$. This means that the line $\ell_1$ through both $u$ and $x$ 
does not lie in $[x, y_1]$  and we find $d(u, \ell) = 2$. Hence $|\Gamma_2(u) \cap \ell| = \phi_2 $. 

As $u \in \Gamma(x)$ and $d(u, y_1) =3$, we find that every vertex in $S_1$ is at distance 2 from $u$. As a consequence we find
$$|\Gamma_2(u) \cap S_2| \leq |\Gamma_2(u) \cap \ell|  - |S_1| =\phi_2-\phi_1.$$
It follows that there are at most $\phi_2-\phi_1$ vertices $v$ in $S_2$ such that $\ell_1 \in [x, v]$, as this implies $d(u, v) \leq 2$. 

On the other hand, for $v \in S_2$, the set $[x, v]$ contains at least $\tau_2-t$ lines through $x$ that are not in $T$. 
By counting pairs $(v, \ell_2)$ with $v \in S_2$ and  lines $\ell_2$ in $[x, v]$ that are not in $T$ we have
$(\beta-\phi_1+1)(\tau_2-t) = n_2(\tau_2-t) \leq (r-t)(\phi_2-\phi_1)$, as $\ell_2 $ contains a vertex $w$ that is in $\Gamma_3(y_1) \cup \Gamma_3(y_2)$. 

Let $f(t) = \frac{r-t}{\tau_2-t} $ for $\phi_1 \leq t \leq \tau_2-1.$ 
As $$f(t) = \frac{r-t}{\tau_2-t}= \frac{r-\tau_2}{\tau_2-t} +1 \leq r-\tau_2+1=f(\tau_2-1),$$ we find 
$$\beta - \phi_1 + 1 \leq \frac{r-t}{\tau_2-t}(\phi_2-\phi_1) \leq (r-\tau_2+1)(\phi_2-\phi_1).$$ 
So, if $\beta > (r-\tau_2+1)(\phi_2-\phi_1) + \phi_1 - 1$, the sets $[x, y_1]$ and $[x, y_2]$ are equal. This proves the theorem.
\qed

The following result was shown in \cite{LK2024+}.

\begin{theorem}[{\cite[Theorem 26]{LK2024+}}]\label{srsg}
Let $\Gamma$ be a geometric distance-regular graph with respect to $\mathcal C$, having valency $k$, smallest eigenvalue $\theta_{\min}$, diameter $D \geq 3$, and $\phi_1 \geq 2$. Let $r = -\theta_{\min}$ and $\beta = \frac{k}{-\theta_{\min}}$.  
Assume that $\Gamma$ satisfies the ELS property. Then the following hold:
Let $x, y$ be a pair of vertices at distance 2 in $\Gamma$. Assume that $[x, y] = \{\ell_1, \ell_2, \ldots, \ell_{\tau_2}\}$. Let $\Sigma= \Sigma(x,y)$ be the subgraph induced on 
$\{x\} \cup \ell_1 \cup \ell_2 \cup \cdots \cup \ell_{\tau_2} \cup \{ z \in \Gamma_2 (x) \mid [x, y]= [x, z]\}$.
Then the graph $\Sigma$ has the following properties:
\begin{enumerate}
\item $\Sigma$ has diameter 2.
\item $\Sigma$ is geodetically closed.
\item If $\Sigma$ contains two vertices of a Delsarte clique $C \in {\mathcal C}$, then $\Sigma$ contains all the vertices of $C$.
\item $\Sigma$ is a geometric strongly regular graph with parameters 
\[
\left( \frac{\beta(\beta - \phi_1 + 1)(\tau_2 - 1)}{\phi_1} + \beta \tau_2 + 1, \ \beta \tau_2, \ \beta - 1 + (\phi_1 - 1)(\tau_2 - 1), \ \phi_1 \tau_2 \right).
\]
\item For every pair of vertices at distance 2, there exists a unique induced subgraph $\Delta$ of $\Gamma$ containing them such that $\Delta$ has diameter 2, is geodetically closed, and has the property that whenever $\Delta$ contains two vertices of a Delsarte clique in $\mathcal{C}$, it  contains the entire Delsarte clique.
\end{enumerate}
\end{theorem}

We end this section with another property of the graphs $\Sigma(x,y)$. 

\begin{lemma}\label{intersectinglines}
Let $\Gamma$ be a geometric distance-regular graph with respect to $\mathcal C$, having  diameter $D \geq 3$, and $\phi_1 \geq 2$. 
Assume that $\Gamma$ satisfies the ELS property. Let $\mathcal{C}^{(2)} = \{ \Sigma(x,y) \mid x, y \in V(\Gamma), \ d(x,y) = 2 \}$, where $\Sigma(x,y)$ is as defined in Theorem~\ref{srsg}.
Then, for each pair of distinct intersecting lines $\ell_1$ and $\ell_2$, there exists a unique element $\Sigma \in \mathcal{C}^{(2)}$ such that $\ell_1 \cup \ell_2 \subseteq V(\Sigma)$.
\end{lemma}
{\bf Proof.}
As $\ell_1$ and $\ell_2$ are both maximum cliques and are distinct, there exists $y_i$ on $\ell_i$ for $i=1,2$ such that $y_1$ and $y_2$ are not adjacent. 
Since $\ell_1$ and $\ell_2$ share a common neighbor of $y_1$ and $y_2$, and $\Sigma(y_1, y_2)$ is geodetically closed by Theorem~\ref{srsg}, it follows that $\ell_1$ and $\ell_2$ are both contained in $\Sigma(y_1, y_2)$. This shows the existence of $\Sigma$.
The uniqueness of $\Sigma$ follows from Theorem~\ref{srsg}~(5).
\qed

\section{Bounds on geometric parameters}
\begin{theorem}\label{2design}
Let $\Gamma$ be a geometric distance-regular graph with respect to $\mathcal C$, having valency $k$, smallest eigenvalue $\theta_{\min}$, diameter $D \geq 3$, and $\phi_1 \geq 2$. Let $r = -\theta_{\min}$.  
If $\Gamma$ satisfies the ELS property, then the following holds:
\begin{enumerate}

\item For each vertex $x\in V(\Gamma)$, let $\mathcal{L}$ be the set of lines through $x$ and $\mathcal{B}=\{[x, z] \mid z\in \Gamma_2(x)\}$. Then $(\mathcal{L}, \mathcal{B})$ is a $2$-$(r, \tau_2, 1)$-design.
\item $r\geq \tau_2(\tau_2-1)+1.$
\item $\phi_2\geq \tau_2(\phi_1-1)+1.$
\end{enumerate}
\end{theorem}
{\bf Proof.}
(1): By Lemma~\ref{Del},  each vertex lies on $r$ lines. Let $ \mathcal{C}^{(2)} = \{ \Sigma(u,v) \mid u, v \in V(\Gamma), \text{ and } d(u,v) = 2 \} $, as defined in Theorem~\ref{srsg}. 
Let $x$ be a vertex, and let $\ell_1$ and $\ell_2$ be two distinct lines that contain $x$. Then by Lemma~\ref{intersectinglines}, there exists a unique element $\Sigma \in \mathcal{C}^{(2)} $ for which $\ell_1 \cup \ell_2 \subseteq V(\Sigma)$. Now, the set $\mathcal{B} = \{ [x, z] \mid z \in \Gamma_2(x) \}$ is in one-to-one correspondence with the elements of $\mathcal{C}^{(2)}$ that contain $x$. This shows that $(\mathcal{L}, \mathcal{B})$ is a $2$-$(r, \tau_2, 1)$-design.

(2):  As $\phi_1 \geq 2$, by Lemma~\ref{tauincrease}, $r=\tau_D > \tau_{D-1} > \cdots > \tau_2 \geq \phi_1 \geq 2$ holds. Now, by Lemma~\ref{subdesign}, we find
 $r\geq \tau_2(\tau_2-1)+1$.


(3):  Let $x,y$  be two vertices at distance 2 in $\Gamma$. Let $\Sigma = \Sigma(x,y)$, as defined in Theorem~\ref{srsg}. Let $\ell_1, \ell_2, \ldots, \ell_{\tau_2}$ be the lines through $x$ 
in $\Sigma$. 
Let $\ell$ be a line through $x$ but not in $\Sigma$. Then $d(y, \ell) = 2$. 
Let $S:=\Gamma(x)\cap \Gamma(y)$ and $T:= (\Gamma_2(y)\cap\ell) \setminus \{x\}$. We will show that 
$|T|\geq \tau_2(\phi_1-1)$ from which it follows that $\phi_2\geq \tau_2(\phi_1-1)+1$ holds.
\\
\indent
Let $\mathcal{L}_y$ be the set of lines through $y$ and $\mathcal{B}_y=\{[y, z]~|~z\in\Gamma_2(y)\}$. By (1), $(\mathcal{L}_y, \mathcal{B}_y)$ is a $2$-$(r, \tau_2, 1)$-design. 
For each vertex $u \in T$, we have $u \notin V(\Sigma)$. Therefore, $[y, u] \neq [y, x]$, and hence $|[y, u] \cap [y, x]| \leq 1$ holds. In particular, $u$ has neighbors on at most one line of $[y, x]$.

By double counting the edges between $S$ and $T$, we have $$|T|\phi_1\geq c_2(\phi_1-1)=\tau_2\phi_1(\phi_1-1).$$
Thus, $|T| \geq \tau_2(\phi_1-1)$ holds. This shows (3).
\qed

\section{Geometric distance-regular graphs with $\phi_1=\tau_2\geq2$}\label{phi=tau}
In this section, we will give a characterization of geometric distance-regular graphs with diameter $D\geq3$ and $\tau_2=\phi_1\geq2$. First we will give an equivalent reformulation of 
\cite[Theorem~4.6]{cuy}. It gives a characterization of certain partial linear spaces and we will use it later to give characterizations of the Grassmann and Johnson graphs. 

\begin{theorem}[cf.~{\cite[Theorem~4.6]{cuy}}]\label{cuyprs-pls}
Let $(\mathcal{P}, \mathcal{L}, \mathcal{I})$ be a connected partial linear space that is not a linear space, and suppose it satisfies the following conditions for some fixed integer $\alpha \geq 2$.

\begin{enumerate}
\item For each pair of nonincident point $p$ and line $\ell$, the point $p$ is collinear with either $0$ or $\alpha$ points lying on $\ell$.

\item For each pair of noncollinear points $p_1$ and $p_2$, we have \\
$\bigl|\{\, p \in \mathcal{P} \mid p \text{ is collinear with both } p_1 \text{ and } p_2 \,\}\bigr| \leq \alpha^2$.

\item There exists a line $\ell$ such that $1+\alpha < |\{\, p \in \mathcal{P} \mid p \text{ is incident with } \ell \,\}| < \infty$.

\item There exists a point that is incident with more than $\alpha$ lines.
\end{enumerate}

Let $q=\alpha-1$. Then one of the following holds.

\begin{enumerate}
\item[(i)] $q$ is a prime power. Let $V$ be a vector space of dimension $d$ (where $d \geq 4$ is a cardinal number) over the field $\mathbb{F}_{q}$; we have $(\mathcal{P}, \mathcal{L}, \mathcal{I}) \cong (W_i, W_{i+1}, \subseteq)$, where $i$ is a finite integer with $1 < i < d-1$, and $W_j$ denotes the set of all $j$-dimensional subspaces of $V$ for $j = i, i+1$.

\item[(ii)] $q=1$. Let $X$ be a set of cardinality $|X|$; we have $(\mathcal{P}, \mathcal{L}, \mathcal{I}) \cong \left( \binom{X}{i}, \binom{X}{i+1}, \subseteq \right)$, where $i$ is a finite integer with $1 < i < |X| - 1$, and $\binom{X}{j}$ denotes the set of all $j$-subsets of $X$ for $j = i, i+1$.
\end{enumerate}

\end{theorem}

\begin{remark}
\begin{enumerate}
\item The dimension $d$ and the cardinality $|X|$ need not be finite.
\item If $d$ is finite, then the incidence structure $(W_i, W_{i+1}, \subseteq)$ is isomorphic to \\ 
$(W_{d-i}, W_{d-i-1}, \supseteq)$ for $0 \leq i \leq d-1$.
\item If $|X|$ is finite, then $\left( \binom{X}{i}, \binom{X}{i+1}, \subseteq \right)$ is isomorphic to $\bigl(\binom{X}{|X|-i}, \binom{X}{|X|-i-1}, \supseteq \bigr)$ for $0 \leq i \leq |X|-1$.
\end{enumerate}
\end{remark}

Now we will give the promised characterization of the Grassmann and Johnson graphs. It is an extension of \cite[Theorem 7.1]{KB2010}, in which we also repair some inaccuracies in the statement of \cite[Theorem 7.1]{KB2010}.

\begin{theorem}\label{JG}
Let $\Gamma$ be a geometric distance-regular graph with diameter $D \geq 3$ and parameters satisfying $\phi_1 = \tau_2 \geq 2$.  Then one of the following holds:
\begin{enumerate}
\item $\phi_1=\tau_2=2$ and $\Gamma$ is a Johnson graph. 
\item $\phi_1=\tau_2\geq3$ and $\Gamma$ is a Grassmann graph defined over the field $\mathbb{F}_{\tau_2-1}$.
\end{enumerate}
\end{theorem}
{\bf Proof.}
Let $\alpha = \phi_1 = \tau_2\geq 2$, and $\mathcal{C}$ be a set of Delsarte cliques of $\Gamma$ such that $\Gamma$ is the point graph of the partial linear space $(V(\Gamma), \mathcal{C}, \in)$. Since $D \geq 3$, there exists a pair of noncollinear vertices, so $(V(\Gamma), \mathcal{C}, \in)$ is not a linear space. 
We will further show that all conditions of Theorem~\ref{cuyprs-pls} are satisfied for $(V(\Gamma), \mathcal{C}, \in)$.

By the definition of the parameter $\phi_1$, Condition~(1) of Theorem~\ref{cuyprs-pls} is satisfied. Moreover, by Equation~\eqref{eqgeom}, we have $c_2 = \tau_2 \phi_1 = \alpha^2$, and therefore Condition~(2) of Theorem~\ref{cuyprs-pls} also holds. By Lemma~\ref{2phi1}, each line $C \in \mathcal{C}$ is incident with at least 
$2\phi_1 \geq \alpha + 2 > \alpha + 1$ vertices. Since $\Gamma$ is finite, Condition~(3) of Theorem~\ref{cuyprs-pls} is satisfied as well. Finally, by Lemma~\ref{tauincrease}, each vertex is incident with  $\tau_D > \tau_2 = \alpha$ lines, so Condition~(4) of Theorem~\ref{cuyprs-pls} is satisfied.
Consequently, all conditions of Theorem~\ref{cuyprs-pls} are fulfilled. Therefore, by Theorem~\ref{cuyprs-pls}, one of the following holds.

\begin{enumerate}
\item[(i)] $\alpha - 1$ is a prime power, and hence $\phi_1 = \tau_2 \geq 3$.  
Let $V$ be a vector space of dimension $d$ (where $d \geq 4$ is a finite integer) over the field $\mathbb{F}_{\alpha - 1}$. Then
\[
(V(\Gamma), \mathcal{C}, \in) \cong (W_i, W_{i+1}, \subseteq),
\]
where $i$ is a finite integer with $1 < i < d - 1$, and $W_j$ denotes the set of all $j$-dimensional subspaces of $V$ for $j = i, i+1$.

\item[(ii)] $\alpha - 1 = 1$, and hence $\phi_1 = \tau_2 = 2$.  
Let $X$ be a finite set of cardinality $|X|$. Then
\[
(V(\Gamma), \mathcal{C}, \in) \cong \bigl(\tbinom{X}{i}, \tbinom{X}{i+1}, \subseteq\bigr),
\]
where $i$ is a finite integer with $1 < i < |X| - 1$, and $\tbinom{X}{j}$ denotes the set of all $j$-subsets of $X$ for $j = i, i+1$.
\end{enumerate}

For the first case, the point graph of $(W_i, W_{i+1}, \subseteq)$ is the graph whose vertices are the $i$-dimensional subspaces of $V$, where two distinct vertices are adjacent if and only if they intersect in an $(i-1)$\nobreakdash-dimensional subspace. Consequently, the point graph is a Grassmann graph defined over the field $\mathbb{F}_{\tau_2 - 1}$.

Similarly, in the second case, the point graph of $\bigl(\binom{X}{i}, \binom{X}{i+1}, \subseteq\bigr)$ is a Johnson graph.

This proves Theorem~\ref{JG}. 
\qed

\bigskip

\section{Dual Pasch axiom}
We begin by recalling the dual Pasch axiom.

Let $\Gamma$ be a graph such that it is the point graph of a partial linear space $X= (V(\Gamma), \mathcal{L}, \in)$. 
Then $\Gamma$ satisfies the \emph{dual Pasch axiom} if for each pair of two adjacent vertices $x, y$ of $\Gamma$, the set of common neighbours of $x$ and $y$, outside the line through 
$x$ and $y$, forms a clique.

The following lemma gives a sufficient condition for a geometric distance-regular graph to satisfy the dual Pasch axiom, and its proof follows \cite{WB1984}.
\begin{lemma}\label{dualpaschlemma}
Let $\Gamma$ be a geometric distance-regular graph with respect to $\mathcal C$, having  diameter $D \geq 3$, and $\phi_1 \geq 2$. If  $\Gamma$ has the ELS property, then $\Gamma$ satisfies the dual Pasch axiom.
\end{lemma}
{\bf Proof.}
As $\Gamma$ has the ELS property, we define $\mathcal{C}^{(2)}$ as we did in Lemma~\ref{intersectinglines}. Then, for each pair of intersecting lines $l_1$ and $l_2$, there exists a unique element $\Sigma \in \mathcal{C}^{(2)}$ such that $l_1 \cup l_2 \subseteq V(\Sigma)$. We denote this element by $\Sigma(l_1, l_2)$. Furthermore, for a pair of adjacent vertices $x$, $y$, we denote the unique line containing both of them by $l_{xy}$.
For a pair of vertices $x$, $y$ at distance 2, following the notation of Lemma~\ref{intersectinglines}, let $\Sigma(x, y)$ be the unique element in $\mathcal{C}^{(2)}$ determined by $x$ and $y$.

Let $x_1$, $x_2$, $y_1$, $y_2$ be four distinct vertices such that
$x_1 \sim x_2$, $y_j \notin l_{x_1 x_2}$, and $x_i \sim y_j$ for $1 \leq i, j \leq 2$. To show that $\Gamma$ satisfies the dual Pasch axiom, we need to show that $y_1 \sim y_2$.

Suppose $y_1 \not\sim y_2$. Both $x_1$ and $x_2$ lie on a geodesic from $y_1$ to $y_2$. By Theorem~\ref{srsg}, each element in $\mathcal{C}^{(2)}$ is geodetically closed. Therefore, $x_1, x_2 \in V(\Sigma(y_1, y_2))$, and hence the line $l_{x_1 x_2}$ is contained in $\Sigma(y_1, y_2)$. This means that all the distinct lines $l_{x_i y_j}$ for $1 \leq i, j \leq 2$ and $l_{x_1 x_2}$ lie in $\Sigma(y_1, y_2)$, and thus $\Sigma(y_1, y_2)$ is the only element of $\mathcal{C}^{(2)}$ that contains any 3 of $x_1$, $x_2$, $y_1$, and $y_2$.

By Lemma~\ref{tauincrease}, we have $\tau_D > \tau_2$, and hence there exists a line $h$ containing $x_1$ such that $h$ intersects $\Sigma(y_1, y_2)$ at the unique vertex $x_1$. Since $\phi_1 \geq 2$, $x_2$ has a neighbor in $h \setminus \{x_1\}$, denoted by $z$. Clearly, $z \notin V(\Sigma(y_1, y_2))$.

As $z \notin V(\Sigma(y_1, y_2))$, without loss of generality, we may assume that $z \not\sim y_1$. Since both $x_1$ and $x_2$ lie on a geodesic between $y_1$ and $z$, we have $x_1, x_2 \in V(\Sigma(y_1, z))$. Therefore, $\Sigma(y_1, z)$ contains both the lines $l_{x_1 y_1}$ and $l_{x_2 y_1}$, and hence $\Sigma(y_1, z) = \Sigma(y_1, y_2)$, a contradiction. 

This shows that $y_1 \sim y_2$, and hence $\Gamma$ satisfies the dual Pasch axiom.
\qed

\vspace{0.2cm}
Let $\Gamma$ be a geometric distance-regular graph with respect to $\mathcal{C}$, satisfying the dual Pasch axiom. We will show some properties of the local graph of such a graph. To do so, we introduce the concept of assemblies. A maximal clique in $\Gamma$ that is not contained in $\mathcal{C}$ is called an \emph{assembly}.

\begin{proposition}\label{assembly}
Let $\Gamma$ be a geometric distance-regular graph with respect to $\mathcal C$, having valency $k$, smallest eigenvalue $\theta_{\min}$, diameter $D \geq 3$, and $\phi_1 \geq 2$. Let $r = -\theta_{\min}$ and $\beta = \frac{k}{-\theta_{\min}}$.
Assume that $\Gamma$ satisfies the dual Pasch axiom.
Then each pair of adjacent vertices $x$ and $y$ lies in a unique assembly, denoted by $M_{xy}$. 
The assembly $M_{xy}$ has order $(\phi_1 - 1)r + 1$, and hence $\beta \geq (\phi_1 - 1)r$. 
Consequently, each vertex lies in exactly $\frac{\beta}{\phi_1 - 1}$ assemblies.
\end{proposition}
{\bf Proof.}
Let $U$ be the set of common neighbours of $x$ and $y$ that are not on the unique line $\ell$ through $x$ and $y$. As $a_1 = \beta - 1 + (\phi_1-1)(r-1)$ and $\ell$ contains exactly 
$\beta +1$ vertices, it follows that $|U| = (\phi_1-1)(r-1)$.
By the dual Pasch axiom, every two distinct vertices in $U$ are adjacent. 
Let $u \in U$. Then $u$ has exactly 
$\phi_1$ neighbours in $\ell$, as $u$ is not on $\ell$. 
If $\phi_1=2$, then the induced subgraph $C$ on $U \cup \{x, y\}$ is a maximal clique of $\Gamma$ as $C$ intersects $\ell$ in exactly $\phi_1=2$ vertices. 
This shows that $x,y$ lie in a unique assembly $M_{xy}$ and that $M_{xy}$ has order 
$(\phi_1-1)(r-1) + \phi_1= (\phi_1-1) r+1$ when $\phi_1=2$.

From now we assume $\phi_1 \geq 3$. 
Let $z$ be a neighbour of $u$ in $l$ distinct from $x$ and $y$, 
and let $v \in U$ be distinct from $u$.
Consider the line $\ell_1$ through $x$ and $u$, and the line $\ell_2$ through $y$ and $u$. Note that $z \notin \ell_1 \cup \ell_2$, as two lines do not intersect in more than one vertex.
The vertex $v$ cannot lie on both lines of $\ell_1$ and $\ell_2$, so we may assume without loss of generality that $v\notin \ell_1$. 
The vertices $z$ and $v$ are common neighbours of $x$ and $u$, both not lying on the line $\ell_1$, and they are distinct since $v$ is not on $\ell$. 
By the dual Pasch axiom, the vertices $v$ and $z$ are adjacent. 
This shows that the subgraph $C$  induced on the set $U \cup \{ z \mbox{ on } \ell \mid z \sim u\}$ is a clique and $C$ is maximal, as $C$ intersects $\ell$ in exactly $\phi_1$ vertices. 
Hence $C$ is the maximal clique 
containing $\{x, y\} \cup U$.  This shows that $x,y$ lie in a unique assembly $M_{xy}$ and that $M_{xy}$ has order 
$(\phi_1-1)(r-1) + \phi_1= (\phi_1-1)r +1$, when $\phi_1 \geq 3$.
So we have shown that each assembly has exactly $(\phi_1-1) r +1$ vertices. This implies that $\beta\geq (\phi_1-1) r$ as a clique has at most $\beta +1$ vertices. 
Therefore each vertex $x$ lies in exactly $\frac{k}{(\phi_1-1)r} = \frac{\beta}{\phi_1-1}$ assemblies, as $k = r\beta$.
This shows the proposition.
\qed

\begin{remark} 
Note that the Grassmann graph $J_q(2D, D)$, where $q$ is a prime power and $D \geq 2$, has $\beta = (\phi_1 - 1)r$, and the assemblies are Delsarte cliques.
\end{remark}

The above proposition shows that through each edge, there is a unique assembly, and that all the assemblies have order $(\phi_1-1) r+1$. 
Now we show that each local graph is a $(\phi_1-1)$-clique extension of a grid.

\begin{theorem}\label{assemthm}
Let $\Gamma$ be a geometric distance-regular graph with respect to $\mathcal C$, having valency $k$, smallest eigenvalue $\theta_{\min}$, diameter $D \geq 3$, and $\phi_1 \geq 2$. Let $r = -\theta_{\min}$ and $\beta = \frac{k}{-\theta_{\min}}$.
For a vertex $x$, let $\Delta_x$ be the local graph at vertex $x$. 
Assume that $\Gamma$ satisfies the dual Pasch axiom.
Let $\ell$ be a line, $M$ an assembly and $u$ be a vertex of $\Gamma$ not in $M$. 
Then the following hold:
\begin{enumerate}
\item If $\ell$ intersects with $M$, then they intersect in exactly $\phi_1$ vertices;
\item If $u$ has a neighbour in $M$, then $u$ has exactly $\phi_1$ neighbours in $M$ and there is exactly one line through $u$ that intersects $M$. 
\end{enumerate}
As a consequence we have that for each vertex $x$ of $\Gamma$, the local graph $\Delta_x$ at vertex $x$ is the $(\phi_1-1)$-clique extension of a $\frac{\beta}{\phi_1-1} \times r$-grid.
\end{theorem}
{\bf Proof.}
Let $M$ be an assembly and $\ell$ be a line in $\mathcal C$ and assume that $M$ and $\ell$ intersect in a vertex $w$. Let $\ell_1:=\ell , \ell_2, \ldots, \ell_r$ be the lines through $w$. 
As $M \setminus \ell_i \neq \emptyset$ for $i=1,2,\ldots, r$ and each vertex outside $\ell_i$ has at most $\phi_1$ neighbours in $\ell_i$, we have $|M \cap \ell_i|\leq \phi_1$ for 
$i=1,2,\ldots,r$. Note that $M \setminus \{w\}= \bigcup\limits_{i=1}^{r}(M \cap \ell_i \setminus \{w\})$. Therefore 
$(\phi_1-1) r=|M \setminus \{w\}|\leq \sum\limits_{i=1}^{r} (|M\cap \ell_i|-1) \leq \sum\limits_{i=1}^{r} (\phi_1-1)= (\phi_1-1) r$, and hence $|M \cap \ell_i| = \phi_1$ for $1\leq i \leq r$. 
In particular, $|M\cap \ell|= \phi_1$. 
This shows Item~(1) of the theorem.

Now let $u$ have a neighbour $v$ in $M$. Let $m$ be the line through $u$ and $v$. Then $|m \cap M| = \phi_1$. Let $w \in M$ such that $w$ is not on $m$. Then $w$ has exactly 
$\phi_1$ neighbours on $m$, and those must be in $M$, since $|m \cap M| = \phi_1 $. This shows that $u$ has exactly $\phi_1$ neighbours in $M$ and that there is a unique line 
through $u$ that intersects $M$. This shows Item~(2) of the theorem. 

Let $x$ be a vertex of $\Gamma$. Let $\ell_1, \ell_2, \ldots, \ell_r$ be the lines through $x$ and $M_1, M_2, \ldots, M_t$, where $t = \frac{\beta}{\phi_1-1}$, be the assemblies through $x$.
Let $1 \leq i, a \leq r$ and $1 \leq j, b \leq t$. Then $|\ell_i \cap M_j| = \phi_1$ and $\ell_i \cap M_j \cap \ell_a \cap M_b = \emptyset$ unless $(i, j) =  (a, b)$. This shows that 
$\Delta_x$ is the $(\phi_1-1)$-clique extension of the $\frac{\beta}{\phi_1-1} \times r$-grid.
\qed

\section{Proof of  Theorem~\ref{main2}}
In this section, we will provide the proof of  Theorem~\ref{main2}. We first prove the special case $2 \leq \phi_1 \leq \tau_2 - 2$, which is stated in Theorem~\ref{main1}.

\begin{theorem}\label{main1}
Let $\Gamma$ be a geometric distance-regular graph with valency $k$, smallest eigenvalue $\theta_{\min}$ and diameter $D \geq 3$. Let $r = -\theta_{\min}$ and $\beta = \frac{k}{-\theta_{\min}}$. If $\phi_1 \geq 2$ and $\phi_1 \neq \tau_2 - 1, \tau_2$, then $\beta < (r - \tau_2+1)(\phi_2 - \phi_1) + \phi_1$.
\end{theorem}
{\bf Proof.}
As $\phi_1 \geq 2$, we have $c_2 = \phi_1 \tau_2 \geq 2$. By Lemma~\ref{tau2geqphi1}, we have $\tau_2 \geq \phi_1$. Since $\phi_1 \neq \tau_2 - 1, \tau_2$, it follows that $\tau_2 - 2 \geq \phi_1 \geq 2$, and hence $\tau_2 \geq 4$.

Assume that $\beta \geq (r - \tau_2+1)(\phi_2 - \phi_1) + \phi_1$. By Theorem~\ref{geometric[]}, $\Gamma$ satisfies the ELS property. Furthermore, by Theorem~\ref{srsg}, there exists an induced subgraph $\Sigma$ of $\Gamma$ which is a geometric strongly regular graph with parameters $(\frac{\beta(\beta -\phi_1 +1)(\tau_2 -1)}{\phi_1} + \beta \tau_2 + 1, \beta \tau_2, \beta-1 + (\phi_1 -1)(\tau_2 -1), \phi_1 \tau_2)$.

Note that the sizes of the Delsarte cliques in $ \Sigma $ and $ \Gamma $ are the same. Let $ {\mathcal C'} $ denote the set of all Delsarte cliques in $ \Sigma $. Then, $ (V(\Sigma), {\mathcal C'}, \in) $ is a partial geometry  $pg(\beta, \tau_2 - 1, \phi_1)$. By Theorem~\ref{SRGbound}, we have 
$\beta \leq (\tau_2 - \phi_1)^2(2\phi_1 - 1)$.

 By Theorem~\ref{2design}, we have $\beta \geq (r - \tau_2+1)(\phi_2 - \phi_1) + \phi_1 > (\tau_2(\tau_2 - 1) + 1 - \tau_2+1 )(\phi_1 - 1)(\tau_2 - 1) > (\tau_2 - 1)^3(\phi_1 - 1) \geq (\tau_2 - \phi_1)^2(3\phi_1 - 3) \geq (\tau_2 - \phi_1)^2(2\phi_1 - 1),$ which leads to a contradiction. Therefore, Theorem~\ref{main1} is proved.
\qed

\bigskip
\noindent
{\bf Proof of Theorem~\ref{main2}.}
Let $\mathcal{C}$ be a set of Delsarte cliques of $\Gamma$, such that $\Gamma$ is the point graph of the partial linear space $X = (V(\Gamma), \mathcal{C}, \in)$.
Since $c_2 \geq 2$, by Lemma~\ref{tau2geqphi1}, we have $1 \leq \phi_1 \leq \tau_2$.

Assume $\phi_1 = 1$. It follows from the definition of geometric distance-regular graphs and Lemma~\ref{Del} that, for each vertex $x \in V(\Gamma)$, the neighborhood $\Gamma(x)$ can be partitioned into a collection of pairwise vertex-disjoint cliques, all of order $\beta$. Since $\phi_1 = 1$, there are no edges between different cliques, and hence $\Delta_x$ is the disjoint union of cliques of order $\beta$.

If $\phi_1 = \tau_2 \geq 2$, then by Theorem~\ref{JG},  $\Gamma$ is either a Johnson graph or a Grassmann graph defined over the field $\mathbb{F}_{\phi_1 - 1}$. The value of $b=\frac{b_1}{\theta_1+1}$ for the Johnson graph and the Grassmann
graph follows from \cite[Theorem~9.1.2]{bcn89} and \cite[Theorem~9.3.3]{bcn89},
respectively.

If $2 \leq \phi_1 \leq \tau_2 - 2$, by Theorem~\ref{main1}, we have the bound $\beta < (r - \tau_2 + 1)(\phi_2 - \phi_1) + \phi_1.$

Finally, we consider the case $\phi_1=\tau_2-1\geq 2$. If $\beta \geq (r - \tau_2 + 1)(\phi_2 - \phi_1) + \phi_1$, then by Theorem~\ref{geometric[]}, $\Gamma$ satisfies the ELS property. 
Furthermore, by Lemma~\ref{dualpaschlemma}, $\Gamma$ satisfies the dual Pasch axiom.
Now, Theorem~\ref{assemthm} shows that for each vertex $x$ of $\Gamma$, the local graph $\Delta_x$ at vertex $x$ is the $(\phi_1 - 1)$-clique extension of a $\frac{\beta}{\phi_1 - 1} \times r$-grid.

Furthermore, assume that $r(\phi_1-1)>2b^2+2b+1$.
Since the local graph $\Delta_x$ is the $(\phi_1-1)$-clique extension of a
$\frac{\beta}{\phi_1-1}\times r$ grid, there exist two cliques $C_1$ and $C_2$
in $\Delta_x$ such that $|C_1|=\beta$, $|C_2|=r(\phi_1-1)$, and
$|C_1\cap C_2|=\phi_1-1$.
By Lemma~\ref{Del}, we have $\beta\geq r(\phi_1-1)$.
By Lemma~\ref{localgraph}, the smallest eigenvalue of $\Delta_x$ satisfies
$\theta_{\min}(\Delta_x)\geq -b-1$.
By \cite[Theorem~3.6]{KPY2011}, we have $(\theta_1+1)(\theta_D+1)<-b_1$,
which implies that $r=-\theta_D>b+1$.
Since $\phi_1\geq 2$, it follows that $b<(r-1)(\phi_1-1)$, and hence
$\phi_1-1<r(\phi_1-1)-b=r(\phi_1-1)-(b+1)+1$.
By Lemma~\ref{clique-int-bound}, we conclude that
$\phi_1-1\leq (b+1)-1$, that is, $\phi_1\leq b+1$.

This proves Theorem~\ref{main2}.
\qed

\bigskip
\begin{flushleft}

{\Large\textbf{Acknowledgments}}\vspace{0.2cm}

J.H. Koolen is partially supported by the National Natural Science Foundation of China (No. 12471335), and the 
Anhui Initiative in Quantum Information Technologies (No. AHY150000).\

\end{flushleft}

\clearpage

\end{document}